\documentclass[12pt,reqno]{amsart}
\usepackage{amssymb}
\usepackage{amsmath}

\usepackage{amssymb, amsmath, amsfonts, latexsym}
\usepackage{enumerate}

\newtheorem{thm}{Theorem}[section]
\newtheorem{cor}[thm]{Corollary}
\newtheorem{lem}[thm]{Lemma}
\newtheorem{prop}[thm]{Proposition}
\newtheorem{defn}[thm]{Definition}
\newtheorem{example}[thm]{Example}
\newtheorem{remarks}[thm]{Remark}

\numberwithin{equation}{section} \theoremstyle{remark}

\title[Exponential convergence in the Wasserstein metric]{Exponential convergence in the Wasserstein metric $W_1$ for one dimensional diffusions}

\author{Lingyan Cheng}
\address{Lingyan Cheng. Institute of Applied Mathematics, Academy of Mathematics and Systems Science,
Chinese Academy of Sciences, 100190, Beijing China. }
\email{chengly@amss.ac.cn}

\author{Ruinan Li}
\address{Ruinan Li. Institute of Applied Mathematics, Academy of Mathematics and Systems Science,
Chinese Academy of Sciences, 100190, Beijing China. }
\email{ruinanli@amss.ac.cn}

\author{Liming Wu}
\address{Liming Wu\\  Laboratoire de Math. CNRS-UMR 6620,
Universit\'e Blaise Pascal, 63177 Aubi\`ere, France. }
\email{Li-Ming.Wu@math.univ-bpclermont.fr}

\newcommand{\dd}{\mathbb{D}}
\newcommand{\ee}{\mathbb{E}}

\newcommand{\nn}{\mathbb{N}}
\newcommand{\rr}{\mathbb{R}}
\newcommand{\pp}{\mathbb{P}}

\def\BB{\mathcal B}

\def\LL{\mathcal L}

\def\e{\varepsilon}

\def\vep{\varepsilon}
\def\<{\langle}
\def\>{\rangle}

\def\beq{\begin{equation}}
\def\nneq{\end{equation}}

\def\bdef{\begin{defn}}
\def\ndef{\end{defn}}

\def\bthm{\begin{thm}}
\def\nthm{\end{thm}}

\def\bprop{\begin{prop}}
\def\nprop{\end{prop}}

\def\brmk{\begin{remarks}}
\def\nrmk{\end{remarks}}

\def\bexa{\begin{example}}
\def\nexa{\end{example}}

\def\blem{\begin{lem}}
\def\nlem{\end{lem}}

\def\bcor{\begin{cor}}
\def\ncor{\end{cor}}

\date{}

\def\bexe{\begin{exe}}
\def\nexe{\end{exe}}

\def\bprf{\begin{proof}}
\def\nprf{\end{proof}}

\def\benu{\begin{enumerate}}
\def\nenu{\end{enumerate}}

\def\bdes{\begin{description}}
\def\ndes{\end{description}}

\begin{document}
\maketitle

\begin{abstract} In this paper, we find some general and efficient sufficient conditions for the exponential convergence $W_{1,d}(P_t(x,\cdot), P_t(y,\cdot) )\le Ke^{-\delta t}d(x,y)$ for the semigroup $(P_t)$ of one-dimensional diffusion. Moreover some sharp estimates of the involved constants $K\ge 1, \delta>0$ are provided. Those general results are illustrated by a series of examples.

\end{abstract}



\section {Introduction}

\subsection{Framework}
Let $I$ be an interval of $\rr$ so that its interior $I^0=(x_0, y_0)$ where $ -\infty \le x_0 < y_0 \le +\infty$. Consider the diffusion generator on $I$ :
$$
\LL  = a(x) \frac{d^2}{d x^2}+ b(x) \frac{d}{dx}$$
with the Neumann boundary condition at $\partial I =\{x_0,y_0\}\cap \rr$, where the coefficients $a,b:I \to \rr$ satisfy  :

{\it
{\bf (H1)} $ 0 < a \in C^1(I)$ and $b$ is Borel measurable and locally bounded on $I$.}

We can write $\LL$ in Feller's form in terms of the scale function $s$ and the speed function $m$ :
$$
\LL =\frac{d}{dm} \frac{d}{ds}
$$
 where $s,m$ are  determined by their derivatives given by
$$
s^\prime(x)=\exp\left( -\int_c^x \frac{b(u)}{a(u)} du\right) \mbox{ and } m^\prime(x) = \frac{1}{a(x)s^\prime(x)}
$$
where $c$ is some fixed point in $I$.

Let $C_0^\infty (I)$ be the space of infinitely differentiable real functions $f$ on $I$ with compact support and $C_{0,N}^\infty(I)$ be the space of all functions $f$ in $C_0^\infty (I)$ such that $ f^\prime |_{\partial I} =0$ (i.e. satisfying the Neumann boundary condition).
Let $(X_t)$ be the diffusion process on the interval $I$ generated by $\LL$ with initial value $x$.

We will assume :

{\it
{\bf (H2)} The diffusion process is non-explosive or equivalently ( \cite{IM})  :
$$
\begin{aligned}
&\int_c^{y_0} s^\prime(x)\left( \int_c^x m^\prime(y) dy\right)dx= +\infty \quad &\mbox{if } y_0 \notin I;\\
&\int_{x_0}^{c} s^\prime(x) \left( \int_x^c m^\prime(y)\ dy\right)dx= +\infty \quad &\mbox{if } x_0 \notin I.
\end{aligned}
$$

{\bf (H3)}  The speed measure  $dm=m^\prime (x)dx$ is finite, i.e. $m(I):=\int_I m'(x)dx<+\infty$.

{\bf (H4)} The generator $\LL$ defined on $ C_{0,N}^\infty(I)$ is essentially self-adjoint on $L^2(I,dm)$ or equivalently (\cite{DJ,EB}) :
$$
s\notin L^2((x_0,c],dm) \mbox{ if } x_0 \notin I;\quad s\notin L^2([c,y_0),dm) \mbox{ if } y_0 \notin I.
$$
}

Throughout this paper we assume that {\bf (H1)}-{ \bf (H4)} are satisfied. In this case, $\mu(dx):= \frac{m^\prime(x)}{m(I)} dx$ is the unique invariant probability measure of the diffusion $(X_t)$. Let $(P_t)$ be the transition semigroup of $(X_t)$, $\LL_2$ be the generator of $(P_t)$ on $L^2(I,\mu)$ with domain $\dd(\LL_2)$, which is an extension of $(\LL, C_{0,N}^\infty(I))$.

\subsection{Exponential convergence in $L^1$-Wasserstein metric}
Given an absolutely continuous and increasing function $\rho: I \to \rr$, let $d_{\rho}(x,y) =|\rho(x)-\rho(y)|$ be the metric on $I$ associated with $\rho$.
Define the Lipschitzian norm of $f: I \to \rr$ with respect to (w.r.t. in short) $d_{\rho}$ as
$$\|f\|_{Lip(\rho)}:=\sup_{x,y \in I, x \neq y} \frac{|f(x)-f(y)|}{|\rho(x)-\rho(y)|}.$$

For any two probability measures $\mu$, $\nu$ on $I$, the Wasserstein distance between $\mu$ and $\nu$ w.r.t. a given metric $d(x,y)$ on $I$ is defined by
$$
W_{1,d} (\mu,\nu)= \inf_{\pi} \iint_{I^2} d(x,y)\pi(dx,dy)
$$
where $\pi $ runs over all couplings of $\mu$ and $\nu$, i.e. all probability measures $\pi$ on $I^2$ with the first and second marginal distributions $\mu$ and $\nu$, respectively. We say that $(P_t)$ is exponential convergent in $W_{1,d_\rho}$ if there exist some constants $K \ge 1$ and $\delta>0$ such that
\beq\label{exp-C}
W_{1,d_\rho} (P_t(x,\cdot),P_t(y,\cdot) ) \le d_\rho(x,y) K e^{-\delta t}, \ \forall x,y\in I,\ t\ge 0.
\nneq
By Kantorovich duality relation, this is equivalent to
\beq\label{exp-C2}
\| P_t\|_{Lip(\rho)}:=\sup_{f:\|f\|_{Lip(\rho)}\le 1}\|P_t f\|_{Lip(\rho)} \le Ke^{-\delta t}, \ t\ge 0.
\nneq

In this paper, we are interested in the exponential convergence of $(P_t)$ in $W_{1,d_\rho}$.  When $\rho\in L^1(I,\mu)$, consider the Banach space
$$C_{Lip(\rho),0}:=\{f: I \to \rr: \|f\|_{Lip(\rho)} < +\infty \mbox{ and } \mu(f)=0\}$$
equipped with the Lipschitzian norm $\|\cdot\|_{Lip(\rho)}$. We have immediately

\bprop
If $(P_t)$ is exponential convergent in $W_{1,d_\rho}$, i.e.  \eqref{exp-C} holds and $\rho\in L^2(I, \mu)$, then the Poisson operator $(-\LL_2)^{-1} =\int_0^{+\infty} P_t dt : C_{Lip(\rho),0}\to C_{Lip(\rho),0}$ is bounded, and
$$
\|(-\LL_2)^{-1} \|_{Lip(\rho)}\le \int_0^\infty \|P_t\|_{Lip(\rho)}dt \le \frac{K}{\delta}<+\infty.
$$
\nprop

Recall the following result in Djellout and the third named author \cite{DW}.
\bthm
Assume {\bf (H1)}-{\bf (H4)} and let $\rho$ be a function on $I$ such that $\rho \in C^1(I) \cap L^2(I,\mu)$, $\rho^\prime >0$ everywhere,
the Poisson operator
$$
(-\LL_2)^{-1}: C_{Lip(\rho),0}(I) \to C_{Lip(\rho),0}(I)
$$
is well defined and bounded if and only if (iff in short)
\beq\label{cW}
c_W(\rho):= \sup_{x \in I } \frac{s^\prime(x)}{\rho^\prime (x)}\int_x^{y_0} [\rho(y) -\mu(\rho)] m^\prime(y) dy < +\infty.
\nneq
In that case, its norm is given by
$$
\|(-\LL_2)^{-1}\|_{Lip(\rho)} =c_W(\rho).
$$
\nthm

In other words a necessary condition for the exponential convergence of $(P_t)$ in $W_{1,d_\rho}$ is $c_W(\rho)<+\infty$.
The objective of this paper is to show that this necessary condition becomes sufficient in a quite general situation and some sharp estimates of the involved constants $\delta, K$ could be obtained.

\subsection{Some comments on the literature}
For the one-dimensional diffusions, the Poincar\'e inequality (equivalent to the exponential convergence in $L^2(\mu)$) and the log-Sobolev inequality (equivalent to the exponential convergence in entropy) can be characterized by means of the generalized Hardy inequality (\cite{BaRo, Chen05}). For the characterization of  Latala-Oleszkiewicz inequality (\cite{LO}) between Poincar\'e and log-Sobolev, see Barthe and Robertho  \cite{BaRo}. See \cite{CW16} for the characterization of  the Sobolev inequality.

For the transport and isoperimetric inequalities, see Djellout and Wu \cite{DW} for sharp estimates of constants.

On the other hand the exponential convergence in the $L^p$-Wasserstein metric is also a very active subject of study. For the early study on this question, the reader is referred to  \cite{Chen92, CW95, CW97, Chen05}.
Renesse and Sturm \cite{RS05} showed that for $\LL=\Delta -\nabla V\cdot\nabla$ on a Riemannian manifold, if the exact exponential convergence (\ref{exp-C})  in $W_{1,d}$ (with $K=1, d$ being the Riemannian metric) holds, then the Bakry-Emery's curvature  must be bounded from below by $\delta$. If one allowed $K>1$  in
 (\ref{exp-C}), Eberle \cite{EB05} found sharp sufficient conditions for  high dimensional  interacting diffusions, by means of reflected coupling.    The reader is referred to this last paper for an overview of literature.

The general idea of proving the exponential convergence in  $W_{1,d}$  for a high dimenioanl diffusion is to use the reflection coupling $(X_t,Y_t)$ and then to compare $d(X_t,Y_t)$ with some one-dimensional diffusion. But curiously a general study on the exponential convergence in $W_{1,d}$ of one-dimensional diffusion is absent : that is the objective of our study.

\subsection{Organization}
Our paper is organized as follows. In the next section, we state the main result and present several corollaries. In Section 3, we provide several examples to illustrate our theorem. In Section 4, we prove the main result in the compact case. The proof of the general case is given in Section 5.

\section{Main Result}
For a function $f$ on $I$, the $\sup$-norm of $f$ is defined by $\|f\|_\infty =\sup_{x \in I} |f(x)|$. Recall that the conditions {\bf (H1)-(H4)} are always assumed.

\subsection{Main theorem}
Our main result in this paper is :
\bthm \label{mr}
Assume $\rho \in C^1(I) \cap L^2(I,\mu)$ and $ \rho^\prime (x) >0,\ \forall x \in I $. Suppose that $c_W(\rho) <+\infty $.
\benu[(a)]
\item If moreover

{\bf (C)} there exists some positive $C^2$-function $\varphi$ such that :
\benu[(1)]
\item for all $x \in I$, $\rho^\prime (x) \le \varphi (x) \le C \rho^\prime (x)$ for some constant $C \ge 1$ ;
\item for some constant $M \ge 0$,
\beq\label{test f}
(a\varphi^\prime +b \varphi)^\prime \le M \rho^\prime
\nneq
in the distribution sense on $I^0 = (x_0,y_0)$,
\nenu
then for any constant $\alpha \in (0,\frac{1}{M})$,
\beq\label{r1}
W_{1,d_\rho} (P_t(x,\cdot),P_t(y,\cdot) ) \le d_\rho(x,y) K e^{-\delta t},\ \forall x,y\in I, \ t\ge 0,
\nneq
where
\beq
\begin{aligned}
&\delta=\delta(\alpha)=\frac{1-\alpha M}{C \alpha +c_W(\rho)}, \ K=K(\alpha)=\frac{1}{\alpha}\left\| \frac{u+\alpha \varphi }{\rho^\prime}\right\|_\infty, \\
&u(x)=s^\prime(x)\int_x^{y_0} [\rho(y) -\mu(\rho)] m^\prime(y) dy.
\end{aligned}
\nneq
\item
For the metric $d_{\tilde{\rho}}$ where
$
\tilde{\rho}^\prime (x) =u(x)
$ given above,
we have
\beq\label{r2}
W_{1,d_{\tilde{\rho}} }(P_t(x,\cdot),P_t(y,\cdot) ) \le \exp\left(-\frac{t}{c_W(\rho)}\right)d_{\tilde{\rho}}(x,y),\ \forall x,y\in I, \ t\ge 0.
\nneq
\nenu
\nthm

\brmk{\rm As $\LL$ is symmetric on $L^2(I,\mu)$,
by the result in \cite{Wu04}, for any metric $d_\rho$, if \eqref{exp-C} holds, the spectral gap $\lambda_1$ of $\LL$ in $L^2(I,\mu)$ satisfies
$$
\delta \le \lambda_1.
$$
In other words, the best possible exponential convergence rate $\delta$ in any Wasserstein metric $W_{1,d\rho}$ of $(P_t)$ is the spectral gap $\lambda_1$.  Chen's variational formula for the spectral gap $\lambda_1$ says (\cite{Chen99})
$$
\lambda_1=\sup_\rho \frac 1 {c_W(\rho)}.
$$
As $\delta=\delta(\alpha)\to\frac 1 {c_W(\rho)}$ when $\alpha\to 0+$, our estimate on the exponential convergence rate $\delta$ is sharp.

If moreover $\lambda_1>0$ is associated with an eigenfunction $\rho$, which can be chosen to be increasing (\cite{Chen92}), then one verifies easily that
$u(x)=\frac 1{\lambda_1}\rho'(x)$. Thus $c_W(\rho)= \frac 1{\lambda_1}$ and we get from Part (b) above that
$$
W_{1,d_{\rho} }(P_t(x,\cdot),P_t(y,\cdot) ) \le \exp\left(-\lambda_1 t\right)d_{\rho}(x,y),\ \forall x,y\in I, \ t\ge 0.
$$
}
\nrmk

\subsection{Corollaries and examples}

We present now a corollary for illustrating the extra condition {\bf (C)} in Part (a) of Theorem \ref{mr}.

\bcor\label{cor22} Assume that $I=\rr $, $a\equiv 1$ and $\rho(x)=x$, i.e. $d_{\rho} $ is the Euclidean metric on $\rr$.  Suppose that
$$
c_W(\rho)=\sup_{x \in \rr} u(x)< +\infty
$$
where $u(x):=e^{V(x)} \int_x^{+\infty} [y -\mu(y)] e^{-V(y)} dy$ and $V$ is some primitive of $-b(x)$. If
the Bakry-Emery's curvature is bounded from below, i.e.
$V^{\prime \prime}\ge - M$ for some constant $M \ge 0$,  then for any $\alpha\in (0,1/M)$,
\beq\label{rr1}
W_{1,d_\rho} (P_t(x,\cdot),P_t(y,\cdot) ) \le |x-y| K e^{-\delta t},\ \forall x,y\in \rr, \ t\ge 0,
\nneq
where
$$
\begin{aligned}
\delta&=\delta(\alpha)=\frac{1-\alpha M}{ \alpha +c_W(\rho)}, \ K=K(\alpha)=\frac{\left\| u+\alpha \right\|_\infty}{\alpha}. \\
\end{aligned}
$$
\ncor
\bprf
It is obtained by Part (a) of Theorem \ref{mr} with $\varphi(x)=1$ simply.
\nprf

\bcor\label{cor23} Let $I=\rr $ and $\rho(x)=x$ (corresponding to the Euclidean metric). Assume $\int_\rr x^2 d\mu(x)<+\infty$. If
 there exist some constants $C_1 \ge 1$ and $L>|\mu(\rho)|$ large enough such that  $B(x)=\int_0^x\frac{b(y)}{a(y)} dy\to -\infty$ as $x\to\pm\infty$ and
\beq\label{bx1}
b(x)\ne0, \ \frac{x-\mu(\rho)}{-b(x)}\le C_1   ,\quad \mbox{if } |x|>L,
\nneq
then $c_W(\rho)<+\infty$ and $(P_t)$ is exponential convergent in $W_{1,d_{\tilde{\rho}}}$, where
$$ \tilde{\rho}'(x)=u(x)=e^{-B(x)}\int_x^{+\infty} (y-\mu(\rho)) \frac 1{a(y)}e^{B(y)} dy .$$
Furthermore, if either  $b'\le M$  for some non-negative constant $M$, or
\beq\label{bx2}
\frac{ x-\mu(\rho)}{-b(x)} \ge \frac{1}{C_1},\quad \mbox{if } |x|>L,
\nneq
then $(P_t)$ is exponential convergent in $W_{1,d_{\rho}}$ with $\rho(x)= x$.
\ncor

\bprf
For $x>L$, by the mean value theorem of Cauchy, there exists some $ \xi \in (x,+\infty)$ such that
$$
\begin{aligned}
u(x)= \frac{-e^{B(\xi)}(\xi-\mu(\rho))}{e^{ B(\xi)} b(\xi)} = \frac{\xi-\mu(\rho)}{-b(\xi)} \le C_1.
\end{aligned}
$$
For $x <-L$, by the similar proof, we also have $u(x) \le C_1$. Then $c_W(\rho) < +\infty$.  Hence by Part (b) of Theorem \ref{mr}, the exponential convergence \eqref{r2} in the metric $W_{1,d_{\tilde\rho}}$ holds.

Furthermore,  if moreover $b'\le M$ for some non-negative constant $M$, applying Theorem \ref{mr}(a) for
$\varphi=1$, we get the exponential convergence in $W_{1,d_{\rho}}$. In the case where (\ref{bx2}) is satisfied,
we have $u(x) \in [1/C_2,C_2]$ for all $x\in\rr$,  for some constant $C_2\ge C_1$ by the mean value theorem of Cauchy. Then
$$\frac 1{C_2}d_{\tilde\rho}\le d_{\rho}\le C_2 d_{\tilde\rho} .$$ Thus we get the exponential convergence in  $W_{1,d_{\rho}}$-metric :
$$
W_{1,d_{\rho}}(P_t(x,\cdot), P_t(y,\cdot))\le C_2^2 \exp\left(-\frac{t}{c_W(\rho)}\right)|x-y|, \ \forall x,y\in\rr, t\ge0.
$$
\nprf

A curious point in Corollary \ref{cor23} above is that our sufficient condition above for the exponential convergence in  $W_1$-metric associated with the Euclidean distance, depends very few upon the volatility coefficient $a(x)$ (except the conditions {\bf (H1)-(H4)}).

We give an example of Corollary \ref{cor23} :
\bexa{\rm
Let $I=\rr $, $\rho(x)=x$, i.e. $d_{\rho} $ is the Euclidean metric on $\rr$.
Consider the generator
$$
\LL  = \frac{d^2 }{dx^2}+b(x)\frac{d}{dx},
$$
where $ b(x) = -V'(x)= -x(2+\sin x)$.

For this example $\mu(dx) =\frac{1}{C} e^{-V(x)}dx$, where $V(x)=x^2-x\cos x + \sin x$, $C$ is the normalization constant . It is easy to see that the hypotheses {\bf (H1)-(H4)} are all satisfied, and
$ \mu (\rho^2)= \frac{1}{C} \int_\rr x^2 \exp(-x^2+x\cos x - \sin x)dx < +\infty $. Notice that the Bakry-Emery's curvature
$$
V^{\prime\prime}(x)= 2+\sin x + x\cos x
$$
is not bounded from below.

Since
$$
B(x)=-V(x)=-x^2+x\cos x - \sin x \to -\infty \mbox{ as } x \to \pm \infty
$$
and when $L>|\mu (\rho)|$ large enough,
$$
\frac{1}{4} \le \frac{ x-\mu(\rho)}{-b(x)}=\frac{ x-\mu(\rho)}{x(2+\sin x)} \le 4, \mbox{ if }|x| \ge L ,
$$
by Corollary \ref{cor23}, $(P_t)$ generated by $\LL$ is exponential convergent in $W_{1,d_{\rho}}$.

}

\nexa

We present now another example to illustrate the extra condition {\bf (C)}.

\bexa
{\rm
Let $I=\rr, a(x)=1,b(x)=-V'(x)$ where $V(x)$ is an even function such that $V'(x)=x^n(n+1+\sin x) $ for
$x \ge 0$, here $n \ge 2 $.

Since
$$
V^{\prime \prime} (x) = n x^{n -1} (n +1 +\sin x)+ x^n \cos x
$$
is unbounded on $[0,+\infty)$, the Bakry-Emery's curvature is unbounded from below.

For $\rho(x) = x$, we see that,
$$
\lim_{x \to +\infty} u(x)=\lim_{x \to +\infty}e^{V(x)} \int_x^{+\infty} (y-\mu(\rho))e^{-V(y)}dy=\lim_{x \to +\infty}\frac{x -\mu(\rho)}{x^n( n +1 + \sin x)}=0.
$$
Similarly $\lim_{x \to -\infty} u(x)=0$, then
$$
c_W(\rho) =\sup_{x \in \rr} u(x) < +\infty.
$$
We choose the following $\varphi$ :
$$
\varphi(x)=
\begin{cases}
\frac{1}{n+1-\sin x}, \quad &\mbox{if } x \le - L;\\
f(x), \quad &\mbox{if } x \in [-L,L];\\
\frac{1}{n+1+\sin x}, \quad &\mbox{if } x \ge L
\end{cases}
$$
where $L>|\mu(\rho)|$ is a positive constant and $ f(x)$ is $C^2[-L,L]$-function such that $\varphi(x)$ is a $C^2$-function on $\rr$.

For $x \ge L$,
$$
\begin{aligned}
(a\varphi' +b\varphi)'&= \varphi^{\prime \prime} -V'\varphi'-V^{\prime \prime} \varphi\\
&=\frac{(n+1) \sin x + \cos^2 x +1}{ (n +1 +\sin x)^3}-n x^{n -1}\\
&\le M_1
\end{aligned}
$$
where $M_1$ is a positive constant.

For $x<-L$, similarly we have $(a\varphi' +b\varphi)' \le M_2 $ for some constant $M_2>0$.

For $ x \in [-L,L]$, by continuity of $(a\varphi' +b\varphi)'$, $(a\varphi' +b\varphi)' \le M_3$ for some constant $M_3>0$.

Then $(a\varphi' +b\varphi)' \le \max \{M_1,M_2,M_3\}$ on $\rr$.
By Theorem \ref{mr}(a),  the exponential convergence \eqref{r1} w.r.t. $d_\rho$ holds.
}
\nexa

\subsection{Main idea}
We explain now the main idea in Theorem \ref{mr}. The crucial point is that in the actual one-dimensional case, we would have formally
the following commutation relation
\beq\label{a11}
(P_tf)' = P_t^D f'
\nneq
where $P_t^D$ is the semigroup generated by $\LL^D g = (a g' +bg)'$. Then by the Kantorovitch duality relation,
the exponential convergence of $(P_t)$ in $W_{1,d_\rho}$ is equivalent to that of $P_t^D$ to $0$ in the Banach space $b_V\BB$ of all Borel-measurable functions $g$ such that
the norm  $\|g\|_V:=\sup_{x\in I} \frac{|g(x)|}{V(x)}<+\infty$, where $V(x)=\rho'(x)$. An easy sufficient condition
to this last exponential convergence is $\LL^ D U \le -\delta U$ for some positive constant $C$ and some function $U$ such that $\rho'\le U\le C\rho'$.

To see the meaning of the necessary condition $c_W(\rho)<+\infty$, notice that $u$ is a particular solution of $\LL^D u =-\rho'$, $u$ should  be bounded by $C \rho'$ if $P_t^D$ converges exponentially to $0$ in the norm $\|\cdot\|_{\rho'}$. Moreover our extra condition {\bf (C)} says simply that there is some function $\varphi\ge \vep \rho'$ in the domain of $\LL^D$ in $b_V\BB$.

However the formal approach above is very difficult to be realized in the general case. It can be realized rigorously in the compact
 case ($I=[x_0,y_0]$) when $a,b$ are quite regular : see \S 4. The general case can be treated by approximation, as the involved constants $\delta, K$ have explicit expressions.

\section{ Several Examples}
Recall that for $\rho(x)=\int_c^x\frac{1}{\sqrt{a(y)}} dy$ the associated metric $d_{\rho}(x,y)=|\rho(x)-\rho(y)|$ is the intrinsic metric $d_X$ of the diffusion $(X_t)$. In this section, we present  several examples and study  their exponential convergence in the $W_1$-metric associated with  the intrinsic distance $d_X$.

\bexa[Ornstein-Uhlenbeck generator]{\rm
Let $I=\rr$. Consider the Ornstein-Uhlenbeck generator
$$
\LL  = \frac{d^2 }{dx^2}- x \frac{d}{dx},
$$
then $\mu(dx)=\frac{1}{2\pi}e^{-\frac{x^2}{2}}dx$ is the invariant probability measure of $\LL$.
For $\rho(x)=x$, we have $u(x)=1$, thus $c_W(\rho) =1$. Since
$ \tilde{\rho}(x)=x=\rho(x)$, by Part (b) of Theorem \ref{mr},  $ \|P_t\|_{Lip(\rho)} \le e^{-t}$ (well known). We know that the spectral gap $\lambda_1$ of $\LL$ is $1$, which shows that Theorem \ref{mr} is sharp.
\
}
\nexa

\bexa \label{exa22}
{\rm
Let $I=\rr $, $b(x)=-V'(x)$ where $V(x)=C_1|x|^r$,  $C_1$ and $r$ are positive constants, and $a\in C^1(\rr)$ which is bounded i.e. $\frac 1C_2 \le  a(x)\le C_2$ for some constant $C_2\ge1$.
Then $\mu (dx) = \frac 1{Za(x)}e^{B(x)} dx$ ($Z$ being the normalization constant), $B(x)=\int_0^x \frac {b(y)}{a(y)} dy$. For $ \rho (x)=x$, we see that $\frac 1{\sqrt {C_2}}d_\rho\le d_X\le \sqrt{C_2}d_\rho$ and
$$
c_W(\rho)= \sup_{x\in \rr}e^{-B(x)}\int_x^{+\infty} \frac1{a(y)}e^{B(y)} (y-\mu(\rho)) dy.
$$
}
\nexa
For this example it is well known that the spectral gap exists (i.e. $\lambda_1>0$) iff $r\ge 1$. About the exponential convergence in $W_1$ associated to the Euclidean metric,  we have the following result :
\bcor In the above Example \ref{exa22}, for $\rho (x) =x$, $(P_t)$ is exponential convergent in $W_{1,d_\rho}$ iff $r\ge 2$.
\ncor
\bprf At first we can check easily that all assumptions {\bf (H1)-(H4)} hold.
For the necessity, we only need to prove that $c_W(\rho) = +\infty$ when $r <2$. By the L'Hospital criterion,
$$
\begin{aligned}
\lim_{x \to +\infty} e^{-B(x)}\int_x^{+\infty} \frac1{a(y)} e^{B(y)} (y-\mu(\rho)) dy
&=\lim_{x \to +\infty}
 \frac{ -e^{B(x)} (x-\mu(\rho))}{-C_1r |x|^{r-1}e^{B(x)}}\\
&=\lim_{x \to +\infty}\frac{1}{ C_1 r (r-1) |x|^{r-2}}= +\infty.
\end{aligned}
$$
By the definition of $c_W(\rho)$, we have $c_W (\rho)=+\infty$.

For the sufficiency, by the proof above we see that $c_W(\rho) < +\infty$ when $r \ge 2 $. Letting $\varphi=1$, we have
$$
(a\varphi' +b\varphi)'=b^\prime=-C_1 r(r-1)|x|^{r-2}\le 0,
$$
then the condition in Part (a) of Theorem \ref{mr} is satisfied with $C=1$ and $M=0$. Hence by Theorem \ref{mr}(a),  $(P_t)$ is exponential convergent in $W_{1,d_\rho}$ i.e. \eqref{r1} holds with $\delta =\frac{1}{c_W(\rho)+\alpha}$ and $K=\frac{c_W(\rho)+\alpha}{\alpha}$, where $\alpha >0$ is arbitrary. By the equivalence between $d_\rho$ and $d_X$, the exponential convergence \eqref{r1} w.r.t. $d_X$ holds with the same $\delta$ above and $K$ replace by $C_2K$.
\nprf

\brmk
{\rm  What happens if $r\in [1,2)$ ? In fact if one takes an odd and increasing  $C^\infty$ -function $\rho$ such that $\rho(x)=e^{c x^{2-r}}$ for $x\ge 1$ where $c>0$ is arbitrary if $r\in (1,2)$  and $0<c<C_1/2$ if $r=1$.

By Cauchy's mean value theorem,
$$
\frac 1C \le \frac{u(x)}{\rho'(x)} \le C
$$
for some constant $C>1$. Therefore by Theorem \ref{mr}(b), when $r\in [1,2)$, $P_t$ is exponentially convergent in $W_{1,d_\rho}$ for $\rho(x)$ given above.
}
\nrmk

\bexa[Jacobi diffusion]{\rm
Let $I=(0,1)$, $a(x)=x(1-x)$ and $b(x)=-x+\frac{1}{2}$, then $\mu(x) =\frac{1}{ \pi\sqrt{x(1-x)}}$. For
$\rho =\frac{\pi}{2}+ \arcsin(2x-1)$, $d_\rho=d_X$ (the intrinsic metric), we  see that
$$
c_W(\rho) =\sup_{x \in (0,1)} \left( \frac{\pi^2}{8} -\frac{1}{2} \arcsin^2 (2x-1) \right) =\frac{\pi^2}{8},
$$
and by calculus we have
$$
(a \rho^{\prime \prime} +b\rho^\prime)^\prime = 0.
$$
We can choose $\varphi=\rho^\prime$, then
$$
(a \varphi^{ \prime} +b\varphi)^\prime \le 0.
$$
Hence the exponential convergence \eqref{r1} w.r.t. $d_X$ holds with
$$
\begin{aligned}
\delta=\frac{8}{ 8\alpha +\pi^2} \mbox{ and }
K=\frac{\pi^2+8\alpha}{8\alpha}
\end{aligned}
$$
where $\alpha >0$ is arbitrary.
}
\nexa

\bexa[Continuous branching process]
{\rm
Let $I=(0,+\infty)$, $a(x)=2 x$ and $b(x)=-2x+1$, then $\mu(x) =\frac{1}{\sqrt{\pi}} \frac{e^{-x}}{\sqrt{x}}$. This process arises as diffusion limit of discrete space branching process. For $ \rho =\sqrt{2x}$, $d_\rho=d_X$, we see that
$$
c_W(\rho) =\sup_{x \in \rr^+} \left( 1-\frac{e^x}{\sqrt{\pi}} \int_x^\infty \frac{e^{-y}}{\sqrt{y}} dy\right) =1,
$$
and by calculus we have
\beq\label{branching1}
(a \rho^{\prime \prime} +b\rho^\prime)^\prime = -\rho^\prime.
\nneq
We note that
$$
\tilde {\rho}^\prime = \rho^\prime,
$$
by Part (b) of Theorem \ref{mr} the exponential convergence
\eqref{r1} w.r.t. $d_X$ holds with
$$
\begin{aligned}
\delta=1 \mbox{ and }
K=1.
\end{aligned}
$$
Moreover from (\ref{branching1}) we see that the increasing function $\rho-\mu(\rho)$ is an eigenfunction of $-\LL$, and its associated eigenvalue  $1$ must be the spectral gap $\lambda_1$. This example shows again that Theorem \ref{mr} is sharp.
}
\nexa

\bexa[Reflected Bessel diffusion process]{\rm
Let $I=(0,1]$, $a(x)=\frac{1}{2}$ and $b(x)=\frac{\beta-1}{2x}$ where $\beta > 1$ is a constant (the dimension), then $\mu(dx) =\beta x^{\beta-1}dx$. For $ \rho =x$, $\sqrt{2}d_\rho= d_X$, we see that
$$
c_W(\rho) =\sup_{x \in (0,1]} \frac{2}{\beta+1} (x-x^2) =\frac{1}{2(\beta+1)} < \frac{1}{4}.
$$
Choose $\varphi=\rho^\prime=1$, then
$$
(a \varphi^{ \prime} +b\varphi)^\prime =b^\prime =- \frac{\beta-1}{2x^2}  \le 0.
$$
Hence by Theorem \ref{mr}, the exponential convergence \eqref{r1} w.r.t. $d_X$ holds with
$$
\begin{aligned}
\delta=\frac{2(\beta+1)}{ 2\alpha (\beta+1)+1} \mbox{ and }
K=\frac{1}{2}\left(\frac{1}{2\alpha(\beta+1)}+1\right)
\end{aligned}
$$
where $\alpha >0$ is arbitrary.}
\nexa

\section{Compact case}
In this section, we prove the main result Theorem \ref{mr} in the compact case i.e. $I=[x_0,y_0]$ is a bounded closed interval of $\rr$.
\subsection{$C^\infty $-case on compact interval}
Assume $a$, $b \in C^\infty [x_0,y_0]$. For all $f \in C_N^\infty [x_0,y_0]$ where $C_N^\infty [x_0,y_0]= \{f\in C^\infty[x_0,y_0] \mbox{ and } f^\prime |_{\{x_0,y_0\}} =0\}$, $ u(t,x) :=P_t f (x) \in C^\infty (\rr^+ \times [x_0,y_0])$ and it satisfies
$$
\begin{cases}
\partial_t u=\LL u=au^{\prime \prime}+b u^\prime;\\
\partial_x u(t,x_0)=\partial_x u(t,y_0) =0 .
\end{cases}
$$
Let $v(t,x)=\partial_x u(t,x)$, it satisfies
\beq\label{PDE}
\begin{cases}
\partial_t v =(av^\prime+bv)^\prime:=\LL^D v;\\
v(t,x_0)= v(t,y_0) =0,
\end{cases}
\nneq
in other words $v(t,x)$ satisfies the Dirichlet boundary condition. For all $g \in C_D^\infty[x_0,y_0] := \{ g\in C^\infty[x_0,y_0] :  g |_{\{x_0,y_0\}} =0\}$, we define $\LL^D$ as follows :
$$
\LL^D g= (ag^\prime+bg)^\prime.
$$
$(\LL^D, C_D^\infty [x_0,y_0]) $ generates a unique $C_0$-semigroup $(P_t^D)$ on the Banach space $C_D[x_0,y_0]$ of continuous functions on $I=[x_0,y_0]$ satisfying $f(x_0)=f(y_0)=0$. Now we show that :
\blem\label{lem41}
For all $f \in C_D[x_0,y_0]$,
\beq\label{semigp}
P_t^D f(x)=\ee^x 1_{[t < \tau_{\partial I}]} f(X_t)D_t
\nneq
where
\benu
\item
$(X_t)$ is the diffusion generated by $\LL$ with the Neumann boundary condition :
\beq\label{Refdiffu}
dX_t=\sqrt{2 a } (X_t) dB_t+b(X_t) dt +dL_t^{x_0}-dL_t^{y_0}
\nneq
where $(B_t)$ is the Brownian motion, $L_t^{x_0}$ (resp. $ L_t^{y_0}$) is the local time of $(X_t)$ at $x_0$ (resp. $y_0$) ;
\item
$$
\tau_{\partial I} =\inf\{ t\ge 0: X_t \in \{x_0,y_0\} \}
$$
is the first hitting time to the boundary ;
\item
for $t < \tau_{\partial I}$,
$$
 D_t =\exp \left(\int_0^t \frac{a^\prime}{\sqrt{2 a}}(X_s) d B_s + \int_0^t b^\prime (X_s) ds - \frac{1}{4}\int_0^t \frac{{a^\prime}^2}{a}(X_s)ds \right).
$$
\nenu
\nlem
\bprf
For $g\in C_D[x_0,y_0]$, let
$$
\tilde{P}_t^D g(x)=\ee^x 1_{[t < \tau_{\partial I}]} g(X_t)D_t=\ee^x  g(X_t^D)D_t
$$
where $ (X^D_t)$ satisfying $X^D_t=X_t$ for $t<\tau_{\partial I}$ and $X^D_t=X_{\tau_{\partial I} }$ for $t\ge\tau_{\partial I}$, is the killed process at $\partial I$.
First, by the Markov property of $(X_t)$, it is easy to see that $(\tilde{P}_t^D)$ is a $C_0$-semigroup on $C_D[x_0,y_0]$.
Then we prove \eqref{semigp}. When $ 0 \le t <\tau_{\partial I}$, for any $g \in C_D^\infty [x_0,y_0]$,
by It\^{o} formula,
$$
\begin{aligned}
d(D_t g(X_t))&=D_t d g(X_t)+g(X_t) d D_t+ d\langle D, g(X)\rangle_t\\
&=D_t \big(a(X_t) g^{\prime \prime}(X_t)dt + b(X_t) g^\prime (X_t)dt\big)\\
&\quad + D_t\big(g(X_t) b^\prime(X_t)dt+ a^\prime (X_t) g^\prime (X_t)dt\big)+ dM_t\\
&=D_t \LL^D g(X_t) dt + dM_t,
\end{aligned}
$$
where
$$
M_t=\int_0^t D_s\left(\sqrt{2a}(X_s)g^\prime(X_s) +\frac{a^\prime}{\sqrt{2a}}(X_s) g(X_s)\right)dB_s
$$ is a martingale. And for $t\ge \tau_{\partial I}$, $D_t g(X^D_t)=0$. Then for $x \in (x_0,y_0)$,
$$
\begin{aligned}
&\quad \tilde{P}^D_tg(x)-g(x)\\
&=\ee^x(D_t g(X_t^D)) -g(x)\\
&=\ee^x \left(\int_0^t 1_{[s < \tau_{\partial I}]} D_s\LL^D g(X_s)ds + M_{t\wedge \tau_{\partial I}}\right)\\
&=\ee^x \int_0^t 1_{[s < \tau_{\partial I}]} D_s \LL^D g(X_s)ds\\
&=\int_0^t \tilde{P}_s^D \LL^D g(x)ds.
\end{aligned}
$$
Then by the uniqueness, $\tilde{P}^D_tg(x)=P^D_tg(x)$.
\nprf

For every $f\in C_N^\infty[x_0,y_0]$, since $ v(t,x)=\partial_x u(t,x)=(P_t f(x))^\prime$ satisfies the partial differential equation \eqref{PDE} with the initial value condition $ v(0,x)=f^\prime (x)$, it is given by
$$
v(t,x) = P_t^D f^\prime(x).
$$
Hence
\beq\label{PPD}
(P_t f)^\prime = P_t^D f^\prime, \mbox{ for all } t\ge0, \ f\in C_N^\infty[x_0,y_0].
\nneq

Recall that for an everywhere positive function $V$, the $V$-norm of $f$ is defined by $\|f\|_{V}=\sup_{x\in I} \frac{|f(x)|}{V(x)}$.

\blem\label{lem13}
If there exists some positive constant $\delta$ and a $C^2$-function $V: [x_0,y_0] \to \rr^+$ such that $C_1 \le \frac{V}{\rho^\prime} \le C_2$ ($C_1,C_2$ are positive constants) and
\beq\label{LDV}
(a V^\prime +bV)^\prime \le -\delta V \mbox{ on } (x_0,y_0),
\nneq
then
\beq\label{Pd}
\|P_t^D \|_V \le e^{-\delta t}, \forall t\ge 0.
\nneq
Moreover, we have
\beq\label{P}
\|P_t \|_{Lip(\rho)} \le K e^{-\delta t}, \forall t \ge 0
\nneq
where $K=\|\frac{V}{\rho^\prime}\|_{\infty}\| \frac{\rho^\prime}{V} \|_{\infty}$.
\nlem

\bprf The only delicate point is that the test function $V$ does not necessarily belong to the domain of definition of the generator $\LL^D$, in fact $V(x_0), V(y_0)$ may be different of $0$.
At first, we prove \eqref{Pd}. For this purpose it is enough to show $P_t^D V(x) \le e^{-\delta t} V(x)$. Let
$$
Y_t= e^{\delta t} D_t V(X_t) I_{[t<\tau_{\partial I }]},
$$
we only need to prove that $(Y_t)$ is a supermartingale. For $ t<\tau_{\partial I }$, by It\^{o} formula and \eqref{LDV},
$$
d Y_t = \delta Y_t dt+ e^{\delta t}D_t (a V^\prime +bV)^\prime(X_t)dt+dM_t \le dM_t
$$
where $M_t$ is a local martingale up to $\tau_{\partial I}$. Then $(Y_t)$ is a supermartingale by Fatou's lemma. Thus by Lemma \ref{lem41},
$$
e^{\delta t} P_t^D V(x)= \ee^x Y_t \le Y_0=V(x)
$$
where \eqref{Pd} follows.

Now we prove \eqref{P},  which is equivalent to
\beq\label{Pbis}
\|P_tf\|_{Lip(\rho)} \le K e^{-\delta t} \|f\|_{Lip(\rho)},\ \forall f\in C^\infty[x_0,y_0].
\nneq
At first  for $f\in C_N^\infty[x_0,y_0]$,
 we have by \eqref{Pd},

$$
\begin{aligned}
\|P_t f\|_{Lip(\rho)}&=\sup_{x\in I} \left|\frac{(P_t f)^\prime}{\rho^\prime} \right|\\
&\le \| P_t^D f^\prime\|_ V \left\| \frac{V}{\rho^\prime} \right\|_\infty\\
&\le e^{-\delta t} \| f^\prime\|_V \left\|\frac{V}{\rho^\prime} \right\|_\infty\\
&\le \left\|\frac{V}{\rho^\prime}\right\|_\infty \left\|\frac{\rho^\prime }{V}\right\|_\infty e^{-\delta t}\|f\|_{Lip(\rho)}\\
&=:K e^{-\delta t}\|f\|_{Lip(\rho)}.
\end{aligned}
$$
Now for every $f\in C^\infty[x_0,y_0]$ and $n \in \nn^+$, let $f_n=f(x_0)+\int_{x_0}^x \psi_n(y) f'(y)dy$ where $\psi_n(x)=1$ for $x\in[x_0+1/n,y_0-1/n]$ and $\psi_n$ is $C^\infty$-smooth, valued in $[0,1]$, with compact support contained in $(x_0,y_0)$. For each $n \in \nn^+$,  as $f_n$ satisfies the Neumann boundary condition and $\|f_n\|_{Lip(\rho)}\le \|f\|_{Lip(\rho)}$, we have
$$|P_tf_n(x)-P_tf_n(y)|\le K e^{-\delta t} \|f\|_{Lip(\rho)}|\rho(x)-\rho(y)|$$
where \eqref{Pbis} follows by letting $n\to\infty$.
\nprf

Now we turn to :

\bprf[Proof of Theorem \ref{mr} in the compact and $C^\infty$-case]
{\bf Part (a).} Since $\left\| \frac{\rho^\prime }{ u+\alpha \varphi}\right\|_\infty\le \frac{1}{\alpha}$,  it is enough to show \eqref{r1} holds with $\delta=\frac{1-\alpha M}{C \alpha +c_W(\rho)}$, $
K=\left\| \frac{u+\alpha \varphi }{\rho^\prime}\right\|_\infty \left\| \frac{\rho^\prime }{ u+\alpha \varphi}\right\|_\infty$.

By Lemma \ref{lem13}, we only need to find a $C^2$-function $V: [x_0,y_0] \to \rr^+$ such that
$C_1\le\frac{V}{\rho^\prime} \le C_2$ for some positive constants $C_1,C_2$ and
\beq\label{LV}
(aV^\prime +b V)^\prime \le -\delta V \mbox{ on } (x_0,y_0).
\nneq
Consider the following equation
\beq\label{Lu}
\LL^D u= -\rho^\prime \mbox{ on } (x_0,y_0),\ u(x_0)=u(y_0)=0.
\nneq
It is explicitly solvable and the unique solution satisfying the Dirichlet boundary condition is
$$
u(x)=s^\prime(x)\left(\int_x^{y_0}  (\rho(y)-\mu(\rho)) m^\prime(y) dy \right).
$$
It is easy to see $u(x)>0$ for all $x \in (x_0,y_0)$.

Since $ (a \varphi^\prime+b \varphi)^\prime \le M \rho^\prime$, for any constant $\alpha \in (0,\frac{1}{M})$, we can choose
$$
V=\alpha \varphi + u.
$$
First notice that $\sup_{x \in I} \frac{u(x)}{\rho^\prime(x)}=c_W(\rho) < +\infty$, we have
$$
0<\alpha \le \frac{V(x)}{\rho^\prime(x)} \le C \alpha+c_W(\rho).
$$
Moreover
$$
\begin{aligned}
(a V^\prime+b V)^\prime \le \alpha M \rho^\prime -\rho^\prime \le -\frac{1-\alpha M}{C \alpha +c_W (\rho)} V =:-\delta V,
\end{aligned}
$$
then by Lemma \ref{lem13}, we get the desired result.

{\bf Part (b).} It is enough to show \eqref{exp-C2} holds with $\delta=\frac{1}{c_W(\rho)}$, $
K=1$.
Consider the unique solution $u$ of \eqref{Lu} with Dirichlet boundary condition.
Since $ c_W(\rho) <+\infty$, then
$$
(a u^\prime+b u^\prime)^\prime= -\rho^\prime\le - \frac{1}{c_W(\rho)} u.
$$
Since $\tilde{\rho}^\prime =u$, by  Lemma \ref{lem13} with $V=u$, we have
\beq\label{P2}
\|P_t^D\|_u \le e^{-\delta t}.
\nneq
Then
$$
\begin{aligned}
\|P_t \|_{Lip(\tilde{\rho})} =\|P_t^D  \|_u \le e^{-\delta t},
\end{aligned}
$$
which is \eqref{r2}. The proof is finished.
\nprf

\subsection{General case on compact interval}
In this subsection, we prove the main result when $a \in C^1[x_0,y_0]$ and $b$ is Borel measurable and bounded on $[x_0,y_0]$.
\bprf[Proof of Theorem \ref{mr}]
{\bf Part (a).}

{\bf First reduction :} $a,b\in C^1[x_0,y_0]$. Taking $a_\e, b_\e \in C^\infty[x_0,y_0]$ such that $a_\e(x) \to a(x)$, $a_\e^\prime(x) \to a^\prime (x)$, $b_\e(x) \to b(x)$, $b_\e^\prime(x) \to b^\prime (x)$ uniformly over $[x_0,y_0]$ as $\e \to 0$ (i.e. $a_\e\to a$ and $b_\e\to b$ in $C^1[x_0,y_0]$).
Since
$$
u_\e(x)=\exp\left( -\int_c^x \frac{b_\e(y)}{a_\e(y)} dy\right)\int_x^{y_0} \frac{1}{a_\e(y)}[\rho(y) -\mu_\e(\rho)] \exp\left( \int_c^y \frac{b_\e(z)}{a_\e(z)} dz\right) dy
$$
($c$ is a fixed constant in $[x_0,y_0]$),
we see that
$u_\e(x) \to u(x)$ uniformly over $[x_0,y_0]$ as $\e \to 0$, then
$ c_W(\rho, \e)$ defined in \eqref{cW} associated with $(a_\e,b_\e)$ converges to $c_W(\rho)$ associated with $(a,b)$ as $\e \to 0$.

Moreover the condition in \eqref{test f} is satisfied for $(a_\e,b_\e)$ with some constant $M_\e>M$ and $ M_\e \to M$ as $\e \to 0$.
By the result of the $C^\infty$-case in Part (a) of Theorem \ref{mr}, the semigroup $(P_t^\e)$ generated by $\LL^\e =a_\e \frac{d^2}{dx^2}+b_\e \frac{d}{dx}$ satisfies
\beq\label{exp-Ce}
\|P_t^\e \|_{Lip(\rho)} \le K_\e e^{-\delta_\e t}
\nneq
where
$
\delta_\e=\delta_\e(\alpha)=\frac{1-\alpha M_\e}{C \alpha +c_W(\rho,\e)} \mbox{ and }
K_\e=K_\e(\alpha)=\left\| \frac{u_\e+\alpha \varphi }{\rho^\prime}\right\|_\infty \left\| \frac{\rho^\prime }{ u_\e+\alpha \varphi}\right\|_\infty,
$
where $\alpha \in (0,\frac{1}{M_\e})$.
Obviously, $\delta_\e \to \delta$ and $ K_\e \to  K$ as $\e \to 0$.

Now we only need to show for all $f \in  C^\infty[x_0,y_0]$,
\beq\label{Pcov}
P_t^\e f(x) \to P_t f(x), \ \forall x \in [x_0,y_0] \mbox{ as }\e \to 0.
\nneq
But the process $(X_t^\e)$ generated by $\LL^\e$ with the refelection Neumann boundary condition converges in law to $(X_t)$ (well-known in the theory of SDE, \cite{IW}).   Then the convergence above holds.

{\bf Second Reduction :} $a\in C^1[x_0,y_0]$, $b$ is Borel measurable and bounded on $[x_0,y_0]$.
It is enough to show that the semigroup $P_t^{(\delta)}$ generated by $\LL$ on $[x_0+\delta, y_0-\delta]$ with the Neumann boundary conditon at $\{ x_0+\delta, y_0 -\delta\}$ satisfies the conclusion, for all $\delta>0$ small enough, since $P_t^{(\delta)}f\to P_tf$ and all involved constants associated with  $P_t^{(\delta)}$ converge to those related to $P_t$ . Thus working on $P_t^{(\delta)}$ if necessary, we may assume without loss of generality that $a(x)$, $b(x)$ are defined on
$[x_0 -\delta ,y_0+\delta]$ and the condition (\ref{test f}) in Part (a) holds on $(x_0-\delta,y_0+\delta)$.

For any $0<\e<\delta$, taking $ p_\e (x) =\frac{1}{\e } p(\frac{x}{\e})$ where $p$ is a positive  $C^\infty$-function on $\rr$ such that its support is contained in $[-1,1]$ and $\int_\rr p(x)dx=1$. Let
$$
b_\e :=\frac{(b \varphi)*p_\e}{\varphi} \in C^1[x_0-\e,y_0+\e]
$$
where $f*g(x)=\int_\rr f(x-y)g(y)dy$ is the convolution.
By condition \eqref{test f}, we have
$$
( a \varphi^\prime + b \varphi )^\prime * p_\e \le M \rho^\prime *p_\e \le (M+\eta_1(\e)) \rho^\prime\ \text{ on }\ [x_0,y_0]
$$
where $\eta_1(\vep)$ is some positive constant tending to $0$ as $\vep\to 0$.
Since  on $I=[x_0,y_0]$,
$$
\begin{aligned}
&(a \varphi^\prime)^\prime * p_\e \ge (a \varphi^\prime)^\prime -\eta_2(\e) \rho^\prime \mbox{ and }(b \varphi)^\prime *p_\e =(b_\e \varphi)^\prime
\end{aligned}
$$
where $\eta_2(\vep)$ is some positive constant tending to $0$ as $\vep\to 0$.
Then
$$
(a \varphi^\prime)^\prime+ (b_\e \varphi)^\prime\le (a \varphi^\prime)^\prime*p_\e +\eta_2(\e) \rho^\prime + (b \varphi)^\prime *p_\e\le (M+\eta(\e)) \rho^\prime
$$
on $[x_0,y_0]$ where $\eta(\e):=\eta_1(\e)+\eta_2(\e)$.
By the result in the first reduction, for the semigroup $(P_t^\e)$ generated by $\LL^\e f = a f^{\prime\prime}+b_\e f'$ with the Neumann boundary condition,  for any $\alpha\in (0,1/M)$ fixed, we have for any $\e>0$ sufficiently small,
 $\alpha<1/(M+\eta(\e))$,
\beq\label{exp-Ce2}
\|P_t^\e \|_{Lip(\rho)} \le K_\e e^{-\delta_\e t}
\nneq
where
$
\delta_\e(\alpha)=\frac{1-\alpha M_\e}{C \alpha +c_W(\rho,\e)} \mbox{ and }
K_\e(\alpha)=\left\| \frac{u_\e+\alpha \varphi }{\rho^\prime}\right\|_\infty \left\| \frac{\rho^\prime }{ u_\e+\alpha \varphi}\right\|_\infty.
$
$$
u_\e(x)=\exp\left( -\int_c^x \frac{b_\e(y)}{a (y)} dy\right)\int_x^{y_0} \frac{1}{a(y)}[\rho(y) -\mu_\e(\rho)] \exp\left( \int_c^y \frac{b_\e(z)}{a(z)} dz\right) dy
$$
($c$ is a fixed constant in $[x_0,y_0]$).

Since $b_\e(x) \to b(x)$ in measure $dx$ on $[x_0,y_0]$ as $\e \to 0$ and $b_\e\ (\e >0)$, $b$ are uniformly bounded,
we see that
$u_\e(x) \to u(x)$ as $\e \to 0$, then
$ c_W(\rho, \e)$ defined in \eqref{cW} associated with $(a,b_\e)$ converges to $c_W(\rho)$ associated with $(a,b)$ as $\e \to 0$. Moreover,
$ \delta_\e  \to \delta $, $K_\e\to K$ as $\e \to 0$.

Now we only need to show for all $f \in  C_N^\infty[x_0,y_0]$,
$$P_t^\e f(x) \to P_t f(x), \ \forall x \in [x_0,y_0] \mbox{ as }\e \to 0.$$
For all $0<\e<\delta $, consider the diffusion $X_t^\e$ on $[x_0,y_0]$ satisfying
$$
d X_t^\e =\sqrt{2 a}(X_t^\e) d B_t+ b_\e(X_t^\e)dt +d L_t^{x_0,\e}- d L_t^{y_0,\e}
$$
with initial value $x$ and the diffusion $X_t^0$ on $[x_0,y_0]$ satisfying
$$
d X_t^0=\sqrt{2 a}(X_t^0) d B_t+ d L_t^{x_0}- d L_t^{y_0}
$$
with the same initial value $x$.
Then for any $f \in C[x_0,y_0]$, we have by Girsanov's formula,
$$
\begin{aligned}
P_t^\e f(x) &= \ee^x\left[ f(X_t^0) \exp\left(  \int_0^t \frac{b_\e}{ \sqrt{2a}}(X_s^0) d B_s- \frac{1}{4} \int_0^t \frac{b_\e^2}{a} (X_s^0) ds\right)\right]\\
&=\ee^x\left[ f(X_t^0) \exp(  Y_t^\e- \frac 12 \<Y^\e\>_t)\right]
\end{aligned}
$$
where $Y_t^\e :=  \int_0^t \frac{b_\e}{ \sqrt{2a}}(X_s^0) d B_s$. Set
$Y_t=\int_0^t \frac{b}{ \sqrt{2a}}(X_s^0) d B_s$.
Since $\{ \<Y^\e\>_t, \e>0 \}$ is uniformly bounded, the family of  exponential martingales
$\{\exp( Y_t^\e-  \frac 12 \<Y^\e\>_t), \e >0 \}$ is uniformly integrable,  then it is enough to show
$$
\lim_{\e \to 0}\ee^x  \<Y^\e- Y\>_t =  \lim_{\e \to 0}\ee^x \int_0^t \frac{(b_\e-b)^2}{2a}(X_s^0)ds =0.
$$
This follows from the dominated convergence theorem,
since $ \{b_\e, \e >0\}$, $b$ are uniformly bounded, $b_\e \to  b$ in measure $dx$,  and the transition probability
$\pp_x(X^0_s\in dx)$ is absolutely continuous for every $s>0$ by the ellipticity assumption {\bf (H1)}.

Taking $\e \to 0$ in both sides of \eqref{exp-Ce2}, we get the desired result.

{\bf Part (b).} The proof is similar to Part (a), we omit the proof.

\nprf

\section{ General non-compact case}
In this section, we consider the general case : $I=(x_0,y_0)$, $I=[x_0,y_0)$ or $I=(x_0,y_0]$. We only prove Theorem
\ref{mr} in the open interval case.

\bprf[Proof of Theorem \ref{mr}] {\bf Part (a).} For $I=(x_0,y_0)$, $\forall n \in \nn^*$, let $ I_n=[x_n, y_n]$, where $x_n<y_n$, $x_n\downarrow x_0$ and $ y_n \uparrow y_0$ as $ n \to +\infty$ so that the point $c$ in the definition of $s'(x)$ and $m'(x)$ belongs to $(x_1,y_1)$. Let $P_t^{(n)}$ be the semigroup generated by
$\LL_n f(x):=\LL f(x)$ for $x\in I_n$ and $f\in C_N^\infty[x_n,y_n]$ (i.e. satisfying the Neunmann  boundary condition at $x_n,y_n$). Let $X_t^{(n)}$ be the diffusion on the interval $I_n$ generated by $\LL_n $ with the Neunmann boundary condition :
$$
dX_t^{(n)}=\sqrt{2 a } (X_t^{(n)}) dB_t+b(X_t^{(n)}) dt +dL_t^{x_n}-dL_t^{y_n}.
$$

By the result in Section 4 for compact case, we have for all $n>0$ and any arbitrary $\alpha \in (0,\frac{1}{M})$,
\beq\label{exp-C3}
\|P_t^{(n)}\|_{Lip(\rho)} \le K_n e^{-\delta_n t}
\nneq
where
$$\begin{aligned}
\delta_n &= \frac{1-\alpha M}{ C \alpha+c_W(\rho,n)},\ c_W(\rho,n)=\sup_{x\in [x_n,y_n]}\frac{u_n(x)}{\rho'(x)},\\
K_n &=\frac{1}{\alpha}\sup_{x \in I_n}\left| \frac{u_n(x)+\alpha \varphi(x)}{ \rho'(x)} \right|,\\
u_n(x)&=s'(x)\int_x^{y_n} [\rho(y) -\mu_n(\rho)] m'(y) dy, \ x\in I_n
\end{aligned}
$$
with $\mu_n(dx)=\frac{1_{I_n}(x) m'(x)}{m(I_n)} dx$. We see that $u_n(x)\to u(x)$ uniformly as $ n \to +\infty$ over the compact interval $I_N$ for any $N\ge 1$ fixed.

Let $\xi_0\in I$ such that $\rho(\xi_0)=\mu(\rho)$. Fix some $N$ and $\delta>0$ so that $[\xi_0-\delta,\xi_0+\delta]\subset I_N$.  Notice that for $n\ge N$,  if $x\in I_n$ and $x\ge \xi_0+\delta$,
$$\aligned
u_n(x) &\le s'(x)\int_x^{y_0} [\rho(y) -\mu_n(\rho)] m'(y) dy\\
&=u(x) +  [\mu(\rho)-\mu_n(\rho)]s'(x)m([x,y_0]).
\endaligned
$$
For the last term above, if $x\ge \xi_0+\delta$

$$\aligned
s'(x)m([x,y_0))&\le \frac 1{\rho(\xi_0+\delta)-\mu(\rho)} s'(x) \int_{x}^{y_0} [\rho(y)-\mu(\rho)] m'(y)dy\\
&=\frac {u(x)}{\rho(\xi_0+\delta)-\mu(\rho)}.
\endaligned$$

Therefore there is some constant $A_1>0$ such that
$s'(x)m([x,y_0))\le A_1 u(x)$ for all $x\ge \xi_0$.

Now if $x\le \xi_0-\delta$,
$$\aligned
u_n(x) &= s'(x)\int_{x_n}^{x} [-\rho(y) +\mu_n(\rho)] m'(y) dy\\
&\le u(x) +[\mu_n(\rho)-\mu(\rho)] s'(x)m((x_0,x])
\endaligned$$
and similarly we have $s'(x)m((x_0,x])\le A_2 u(x)$ for all $x\le \xi_0$ and some constant $A_2>0$.  Summarizing the previous discussion we get
$$
u_n(x)\le u(x) + |\mu(\rho)-\mu_n(\rho)| Au(x), \ x\in I_n
$$
where $A =\max\{A_1, A_2\}$,
which implies that
$$\limsup_{n\to\infty}c_W(\rho,n)\le c_W(\rho),\ \limsup_{n\to\infty}K_n\le K.$$
(Notice that $\liminf_{n\to\infty}c_W(\rho,n)\ge c_W(\rho),\ \liminf_{n\to\infty}K_n\ge K$ hold always.)

Now for the exponential convergence in Part (a), it remains to show that for any $f \in C_0^\infty (x_0,y_0)$ and $x\in (x_0,y_0)$,
$$
\lim_{n\to \infty } P_t^{(n)}f(x) = P_tf(x).
$$
Denote the first hitting time of $(X_t)$ to the boundary of $I_n$ by
$$
\tau_{\partial I_n}= \inf \{ t\ge 0: X_t \in \{x_n,y_n\} \},
$$
we have
$$
X_t^{(n)}=X_t, \ \forall t \in [0,\tau_{\partial I_n}).
$$
By the non-explosion assumption {\bf (H2)}, we have for any $t\ge0$ and $x\in I$ fixed,
\beq\label{limtau}
\lim_{n \to \infty} \pp_x (\tau_{\partial I_n} >t)=1,
\nneq
then
$$
\begin{aligned}
| P_t^{(n)}f(x)-P_t f(x)|&=|\ee^x f(X_t^{(n)})-\ee^x f(X_t) |\\
&\le 2\|f\|_\infty\pp_x(t \ge \tau_{\partial I_n}) \to 0 \mbox{ as } n \to \infty,
\end{aligned}
$$
by \eqref{limtau}.

\medskip
{\bf Part(b).}  For  any $f \in C_0^\infty (x_0,y_0)$ and $x<y$ in $(x_0,y_0)$, the support of $f$ and $\{x,y\}$ are contained in $I_N$ for some $N\ge1$ large enough. Then if $N\le n\to \infty$, recalling that $\tilde \rho'(x)=u(x)$,
$$\sup_{x\in (x_n,y_n)}\frac{|f'(x)|}{u_n(x)}= \sup_{x\in I_N} \frac{|f'(x)|}{u_n(x)}\to  \sup_{x\in I_N} \frac{|f'(x)|}{u(x)} =\|f\|_{Lip(\tilde \rho)}$$
for $u_n\to u$ uniformly over $I_N$. Since $\displaystyle \lim_{n\to+\infty} c_W(\rho,n)= c_W(\rho)$ as shown in Part (a), we get by Part (b) in the compact case,
$$\aligned
|P_tf(x)-P_tf(y)|&=\lim_{n\to\infty}|P_t^{(n)}f(x)- P_t^{(n)}f(y)|\\
&\le \lim_{n\to\infty} \exp\left(-\frac {t}{c_W(\rho,n)}\right) \int_{x}^y u_n(r) dr  \sup_{z\in (x_n,y_n)}\frac{|f'(z)|}{u_n(z)}\\
&= \exp\left(-\frac {t}{c_W(\rho)}\right)[\tilde\rho(y)-\tilde \rho(x)] \|f\|_{Lip(\tilde \rho)}
\endaligned$$
where the desired result follows by Kantorovitch duality characterization.

\nprf

\end{document}